\documentclass[12pt]{amsart}
\usepackage{amsfonts}
\usepackage{amsmath}
\usepackage{amssymb}
\usepackage[margin=1.2in]{geometry}
\usepackage{bbm}
\usepackage{tikz}
\allowdisplaybreaks[4]

\DeclareMathOperator{\supp}{supp}

\newtheorem{theorem}{Theorem}[section]

\newtheorem{lemma}[theorem]{Lemma}

\theoremstyle{definition}
\newtheorem{definition}[theorem]{Definition}

\newtheorem{remark}[theorem]{Remark}
\newtheorem{Conjecture}[theorem]{Conjecture}

\numberwithin{equation}{section}
\numberwithin{theorem}{section}

\numberwithin{equation}{section}

\begin{document}
\title[Small gaps between almost-twin primes]{Small gaps between almost-twin primes}

\author[Bin Chen]{Bin Chen}
\address{B. Chen\\ Department of Mathematics: Analysis, Logic and Discrete Mathematics\\ Ghent University\\ Krijgslaan 281\\ B 9000 Ghent\\ Belgium}
\email{bin.chen@UGent.be}

\begin{abstract} Let $m \in \mathbb{N}$ be large. We show that there exist infinitely many primes $q_{1}< \cdot\cdot\cdot < q_{m+1}$ such that
	\[
	q_{m+1}-q_{1}=O(e^{7.63m})
	\]
	and $q_{j}+2$ has at most
	\[
	\frac{7.36m}{\log 2} + \frac{4\log m}{\log 2} + 21
	\]
prime factors for each $1 \leq j \leq m+1$.
This improves the previous result of Li and Pan, replacing  $m^{4}e^{8m}$ by $e^{7.63m}$ and $\frac{16m}{\log 2} + \frac{5\log m}{\log 2} + 37$ by $\frac{7.36m}{\log 2} + \frac{4\log m}{\log 2} + 21$. The main inputs are the Maynard-Tao sieve, a minorant for the indicator function of the primes constructed by Baker and Irving, for which a stronger equidistribution theorem in arithmetic progressions to smooth moduli is applicable, and Tao's approach previously used to estimate $\sum_{x \leq n < 2x} \mathbf{1}_{\mathbb{P}}(n)\mathbf{1}_{\mathbb{P}}(n+12)\omega_{n}$, where $\mathbf{1}_{\mathbb{P}}$ stands for the characteristic function of the primes and $\omega_{n}$ are multidimensional sieve weights.
\end{abstract}

\maketitle

\section{Introduction} \label{Intro}
Let $k \in \mathbb{N}$. We consider a set $\mathcal{H}=\{h_{1},...,h_{k}\}$ of distinct nonnegative integers. We call such a set admissible if, for every prime $p$, the number of distinct residue classes modulo $p$ occupied by $h_{i}$ is less than $p$. The following conjecture is one of the greatest open problems in prime number theory.

\begin{Conjecture} [Prime $k$-tuples conjecture] \label{KTC} Given an admissible set $\mathcal{H}=\{h_{1},...,h_{k}\}$, there are infinitely many integers $n$ for which all $n+ h_{i}$ are prime.
\end{Conjecture}

Work on approximations to this conjecture has been very successful in establishing the existence of small gaps between primes. For any natural number $m$, let $H_{m}$ denote the quantity
$$
H_{m}:=\liminf_{n \rightarrow \infty} (p_{n+m} - p_{n}),
$$
where $p_{n}$ denotes the $n$-th prime. In 2013, Zhang \cite{Z2014} proved
$$
H_{1} < 7 \times 10^{7},
$$
by refining the GPY method \cite{GPY} and employing a stronger version of the Bombieri-Vinogradov theorem that is applicable when the moduli are smooth numbers. After Zhang's breakthrough, a new higher rank version of the Selberg sieve was developed by Maynard \cite{J2015} and Tao. This new sieve method provided an alternative way of proving small gaps between primes and had additional consequences. It was more flexible and could show the existence of clumps of primes in intervals of bounded length. Specifically, utilizing the Maynard-Tao sieve, one can show (see \cite{J2015})
$$
H_{m} \ll m^{3}e^{4m},
$$
for any $m\geq 1$. The bound $m^{3}e^{4m}$ was improved by Polymath \cite[Theorem 4(vi)]{PM2014} to $me^{(4-\frac{28}{157})m}$ by incorporating Zhang's version of the Bombieri-Vinogradov theorem \cite{PM2014a}.
In 2017, Baker and Irving \cite{BI2017} achieved a further improvement in the bound:
$$H_{m} \ll e^{3.815m},$$ 
accomplished by constructing a minorant for the indicator function of the primes (see Lemma \ref{BIM} below). This minorant is associated with a stronger equidistribution theorem in arithmetic progressions with smooth moduli. 

One can understand the strengths and the limitations of the current sieve methods by establishing conditional results about primes gaps. Under the Elliott-Halberstam conjecture \cite{EH}, Maynard \cite{J2015} obtained the bound 
$$
H_{1} \leq 12,
$$
improving upon the previous bound $H_{1} \leq 16$ of Goldston, Pintz, and Yıldırım \cite{GPY}.  Let $\mathcal{H}=\{h_{1},...,h_{5}\}=\{0, 2, 6, 8, 12\}$. The proof of the above result relies on considering the quantity
\[S=\sum_{\substack{x \leq n < 2x \\ n\equiv v \, (\! \bmod  W) }} \left( \sum_{j=1}^{5} \mathbf{1}_{\mathbb{P}}(n+h_{j})-1\right)\left( \sum_{\substack{d_{i} | n+h_{i}, 1\leq i \leq 5	 }} \lambda_{d_{1},\cdot\cdot\cdot, d_{5}}\right)^{2},
\]
where $v$ and $W$ are some integers depending on $x$, and 
$$\lambda_{d_{1},\cdot\cdot\cdot, d_{5}}=\mu(d_{1})\cdot\cdot\cdot \mu(d_{5})f\left(\frac{\log d_{1}}{\log R}, \cdot\cdot\cdot ,\frac{\log d_{5}}{\log R}\right)$$ 
with $R=x^{1/2-\varpi}$ for some small $\varpi>0$. The smooth function $f\!: [0,\infty)^{5} \rightarrow \mathbb{R}$ is supported on the simplex\footnote{Throughout the paper, for $k \in \mathbb{N}^{+}$ and $y >0$, we denote by $\Delta_{k}(y)$ the simplex  $\{(t_{1},\cdot\cdot\cdot,t_{k})\in [0,\infty)^{k}:t_{1}+\cdot\cdot\cdot +t_{k}\leq y\}$.}$\Delta_{5}(1)$, and $\mu$ denotes the M$\ddot{\mbox{o}}$bius function. Since the multidimensional sieve weights $\left( \sum \lambda_{d_{1},\cdot\cdot\cdot, d_{5}}\right)^{2} \geq 0$, the inequality $S>0$ would imply that there is some $n \in [x, 2x)$ for which at least two of $n+h_{1}, \cdot\cdot\cdot, n+h_{5}$ are simultaneously prime, and hence there are $2$ primes contained in an interval of length $h_{5}-h_{1}=12$. Maynard's bound $H_{1} \leq 12$  follows readily from establishing $S>0$ for sufficiently large $x$. 

One might anticipate improving the value $12$ to $10$ by delving further into the following sum
\begin{align*}
	S'&=\sum_{\substack{x \leq n < 2x \\ n\equiv v \, (\! \bmod  W) }} \left( \sum_{j=1}^{5} \mathbf{1}_{\mathbb{P}}(n+h_{j})-\mathbf{1}_{\mathbb{P}}(n)\mathbf{1}_{\mathbb{P}}(n+12)-1\right)\left( \sum_{\substack{d_{i} | n+h_{i}, 1\leq i \leq 5	 }} \lambda_{d_{1},\cdot\cdot\cdot, d_{5}}\right)^{2}
	\\
	&=S-\sum_{\substack{x \leq n < 2x \\ n\equiv v \, (\! \bmod  W) }} \mathbf{1}_{\mathbb{P}}(n)\mathbf{1}_{\mathbb{P}}(n+12)\left( \sum_{\substack{d_{i} | n+h_{i}, 1\leq i \leq 5	 }} \lambda_{d_{1},\cdot\cdot\cdot, d_{5}}\right)^{2}.
\end{align*}
To show $S'>0$, it is essential to derive an estimate for $\sum_{n} \mathbf{1}_{\mathbb{P}}(n)\mathbf{1}_{\mathbb{P}}(n+12)\left( \sum \lambda_{d_{1},\cdot\cdot\cdot, d_{5}}\right)^{2}$. Although obtaining an asymptotic estimate for this sum seems to be out of reach by present methods, Tao has devised a method (unpublished) to establish an upper bound. Specifically, the key observation is that, when $n$ and $n+12$ are both primes, we have $ \sum_{\substack{d_{i} | n+h_{i}, 1\leq i \leq 5	 }} \lambda_{d_{1},\cdot\cdot\cdot, d_{5}}= \sum_{\substack{d_{i} | n+h_{i}, 1\leq i \leq 5	 }} \tilde{\lambda}_{d_{1},\cdot\cdot\cdot, d_{5}},$ whenever $\tilde{\lambda}_{d_{1},\cdot\cdot\cdot, d_{5}}=\mu(d_{1})\cdot\cdot\cdot \mu(d_{5})\tilde{f}\left(\frac{\log d_{1}}{\log R}, \cdot\cdot\cdot ,\frac{\log d_{5}}{\log R}\right)$ for another real-valued smooth function $\tilde{f}$ that satisfies $\supp \tilde{f} \subseteq \Delta_{5}(1)$ and $\tilde{f}(0,t_{2}, t_{3}, t_{4}, 0)=f(0,t_{2}, t_{3}, t_{4}, 0)$.
We therefore have 
\begin{align*}
	&\sum_{\substack{x \leq n < 2x \\ n\equiv v \, (\! \bmod  W) }} \mathbf{1}_{\mathbb{P}}(n)\mathbf{1}_{\mathbb{P}}(n+12)\left( \sum_{\substack{d_{i} | n+h_{i}, 1\leq i \leq 5	 }} \lambda_{d_{1},\cdot\cdot\cdot, d_{5}}\right)^{2}
	\\
	=&\sum_{\substack{x \leq n < 2x \\ n\equiv v \, (\! \bmod  W) }} \mathbf{1}_{\mathbb{P}}(n)\mathbf{1}_{\mathbb{P}}(n+12)\left( \sum_{\substack{d_{i} | n+h_{i}, 1\leq i \leq 5	 }} \tilde{\lambda}_{d_{1},\cdot\cdot\cdot, d_{5}}\right)^{2}
	\\
	\leq& \sum_{\substack{x \leq n < 2x \\ n\equiv v \, (\! \bmod  W) }} \mathbf{1}_{\mathbb{P}}(n)\left( \sum_{\substack{d_{i} | n+h_{i}, 1\leq i \leq 5	 }} \tilde{\lambda}_{d_{1},\cdot\cdot\cdot, d_{5}}\right)^{2}.
\end{align*}
According to the Maynard-Tao sieve, the last expression is $$\frac{(1+o(1))x}{\log x(\log R)^{4}}\frac{W^{4}}{\phi(W)^{5}}\int_{\Delta_{4}(1)} \left(\frac{\partial^{4}\tilde{f}(0, t_{2},t_{3}, t_{4},t_{5})}{\partial t_{2}\partial t_{3}\partial t_{4}\partial t_{5}} \right)^{2}\:\mathrm{d} t_{2}\:\mathrm{d} t_{3}\:\mathrm{d} t_{4}\:\mathrm{d}t_{5},$$
as $x \rightarrow \infty$, where $\phi$ represents the Euler totient function. We then need to optimize the square-integral of $\frac{\partial^{4}\tilde{f}(0, t_{2},t_{3}, t_{4},t_{5})}{\partial t_{2}\partial t_{3}\partial t_{4}\partial t_{5}}$ subject to $\tilde{f}$ being supported on $\Delta_{5}(1)$ and having the same trace as $f$ on the boundary $t_{1}=t_{5}=0$. By a converse to Cauchy-Schwarz and the fundamental theorem of calculus, we find that  
\begin{align*}
	&\int_{\Delta_{4}(1)} \left(\frac{\partial^{4}\tilde{f}(0, t_{2},t_{3}, t_{4},t_{5})}{\partial t_{2}\partial t_{3}\partial t_{4}\partial t_{5}} \right)^{2}\:\mathrm{d} t_{2}\:\mathrm{d} t_{3}\:\mathrm{d} t_{4}\:\mathrm{d}t_{5}
	\\
	=&\int_{\Delta_{3}(1)}\:\mathrm{d} t_{2}\:\mathrm{d} t_{3}\:\mathrm{d} t_{4}\int_{0}^{1-t_{2}-t_{3}-t_{4}} \left(\frac{\partial^{4}\tilde{f}(0, t_{2},t_{3}, t_{4},t_{5})}{\partial t_{2}\partial t_{3}\partial t_{4}\partial t_{5}} \right)^{2}\:\mathrm{d}t_{5}
	\\
	\geq&\int_{\Delta_{3}(1)}\:\mathrm{d} t_{2}\:\mathrm{d} t_{3}\:\mathrm{d} t_{4}\frac{1}{1-t_{2}-t_{3}-t_{4}}\left(\int_{0}^{1-t_{2}-t_{3}-t_{4}} \frac{\partial^{4}\tilde{f}(0, t_{2},t_{3}, t_{4},t_{5})}{\partial t_{2}\partial t_{3}\partial t_{4}\partial t_{5}} \:\mathrm{d}t_{5}\right)^{2}
	\\
	=&\int_{\Delta_{3}(1)}\frac{1}{1-t_{2}-t_{3}-t_{4}}\left( \frac{\partial^{3}f(0, t_{2},t_{3}, t_{4},0)}{\partial t_{2}\partial t_{3}\partial t_{4}} \right)^{2}\:\mathrm{d} t_{2}\:\mathrm{d} t_{3}\:\mathrm{d} t_{4}.
\end{align*}
It is clear that the extremizer occurs with
\[\frac{\partial^{4}\tilde{f}(0, t_{2},t_{3}, t_{4},t_{5})}{\partial t_{2}\partial t_{3}\partial t_{4}\partial t_{5}}=\frac{-1}{1-t_{2}-t_{3}-t_{4}}\frac{\partial^{3}f(0, t_{2},t_{3}, t_{4},0)}{\partial t_{2}\partial t_{3}\partial t_{4}}.
\]
While some numerical calculations show that the extremal bound is not strong enough to derive $S'>0$ for all large $x$ (and consequently $H_{1} \leq 10$), Tao's approach remains of independent interest and proves valuable in various applications. To exemplify its significance, we now explore a more general expression than $S$, denoted as $\mathcal{T}$:
\[ \sum_{\substack{x \leq n < 2x \\ n\equiv v \, (\! \bmod  W) }} \mathbf{1}_{\mathbb{P}}(n+h_{j})\mathbf{1}_{\mathbb{P}}(n+h_{\ell})\left( \sum_{\substack{d_{i} | n+h_{i}, 1\leq i \leq k	 }}\mu(d_{1})\cdot\cdot\cdot \mu(d_{k})f\left(\frac{\log d_{1}}{\log x}, \cdot\cdot\cdot ,\frac{\log d_{k}}{\log x}\right)\right)^{2}, 
\]
where the symmetric smooth function $f\!: [0,\infty)^{k} \rightarrow \mathbb{R}$ is supported on the simplex $\Delta_{k}(\delta)$ for some small $\delta > 0$, and $h_{j}, h_{\ell}$ are two distinct element belonging to some admissible set $\mathcal{H}=\{h_{1},...,h_{k}\}$. Applying Tao's approach (repeating the above argument but restricting $\tilde{f}$ to be supported on $\Delta_{k}(\frac{1}{4}-\delta)$ ) and the Maynard-Tao sieve along with the Bombieri-Vinogradov theorem, one can obtain\footnote{For a slightly modified proof, refer to \cite[p. 528, Proof of Lemma 4.6 (iii)]{BFJ2016} }
$$\mathcal{T}\leq (4+O(\delta))\frac{x}{(\log x)^{k}}\frac{W^{k-1}}{\phi(W)^{k}}\int_{\Delta_{k-2}(\delta)} \left(\frac{\partial^{k-2}f(t_{1},\cdot\cdot\cdot t_{k-2},0, 0)}{\partial t_{1}\cdot\cdot\cdot\partial t_{k-2}} \right)^{2}\:\mathrm{d} t_{1}\cdot\cdot\cdot\:\mathrm{d} t_{k-2}.$$
Establishing an upper bound for $\mathcal{T}$ plays a crucial role in exploring the limit points of normalized primes gaps. For a detailed discussion in this direction, we refer to \cite{BFJ2016, JP2018, JM2019}.

In fact, employing the Maynard-Tao sieve, one can derive a small-gaps type result for any subsequence of primes which satisfies the Bombieri-Vinogradov type mean value theorem. Let 
\[ \mathcal{P}_{d}^{(2)} = \{ p: p \ \mbox{is prime and}\  \Omega(p+2) \leq d\},
\]
where $\Omega(n)$ denotes the number of prime divisors of $n$ counted with multiplicity. While there is currently no Bombieri-Vinogradov type mean value theorem (or even asymptotic formula\footnote{  Chen's celebrated theorem \cite{Chen-1978} asserts that $|\mathcal{P}_{2}^{(2)} \cap [1,x]| \gg \frac{x}{(\log x)^{2}}$.}) known for $\mathcal{P}_{d}^{(2)}$,  Li and Pan \cite{L-P2015} successfully established in 2015 small gaps  between primes in $\mathcal{P}_{d}^{(2)}$ by combining ideas from \cite{HB1997, H-T2006, J2015PAP} and the Maynard-Tao sieve. The main objective of the paper is to improve Li and Pan's result by proving the following theorem.

\begin{theorem} \label{NewTh}Let $m \in \mathbb{N}$ be large. Then there exist infinitely many primes $q_{1}< \cdot\cdot\cdot < q_{m+1}$ such that
	\[
	q_{m+1}-q_{1}=O(e^{7.63m})
	\]
	and $q_{j}+2$ has at most
	\[
	\frac{7.36m}{\log 2} + \frac{4\log m}{\log 2} + 21
	\]
	prime divisors for each $1 \leq j \leq m+1$.
\end{theorem}

Li and Pan previously obtained estimates $q_{m+1}-q_{1}=O(m^{4}e^{8m})$ and $\Omega(q_{j}+2)\leq \frac{16m}{\log 2} + \frac{5\log m}{\log 2} + 37$ for every $m \geq 1$. Our improvement on the primes gaps $q_{m+1}-q_{1}$ is based on Baker and Irving's minorant for the indicator function of the primes. Furthermore, the incorporation of Tao's approach and a more meticulous analysis of Li and Pan's argument lead to the sharper estimation of the number of the prime divisors. More specifically, we will utilize Tao's approach to investigate the expression$$ \sum_{\substack{x \leq n < 2x \\ n\equiv v \, (\! \bmod  W) }}  \mathbf{1}_{\mathbb{P}}(n+h_{1})\tau(n+h_{2})\left( \sum_{\substack{d_{i} | n+h_{i}, 1\leq i \leq 2k_{0}	 }} \lambda_{d_{1}, \cdot\cdot\cdot, d_{2k_{0}}}\right)^{2}$$ for large $k_{0}$ (see Section \ref{AUBS2} below for details), where $\tau$ denotes the divisor function.

In this paper, we represent the $2k_{0}$-tuple of real numbers $(x_{1}, \cdot\cdot\cdot, x_{2k_{0}})$ as $\underline{x}$. The greatest common divisor of integers $a$ and $b$ is denoted as $(a,b)$. Additionally, the least common multiple of integers $a$ and $b$ is denoted as $[a,b]$.

\section{Lemmas} \label{LMS}

In this section we introduce two prerequisite results which are quoted from the literature directly. These lemmas play important roles in the proof of our main theorem in Section \ref{PMTHm}.

Our first lemma allows us to choose an admissible set with small gaps, that potentially gives twin primes.  It is a simple application of the Jurkat-Richert theorem \cite[Theorem 8.4]{HR1974} in the theory of sieves.

\begin{lemma} [H. Li and H. Pan {\cite[Lemma 3.1]{L-P2015}}] \label{adsetl} For $k_{0}\geq 1$, there exist $h_{1} < h_{2} < \cdot\cdot\cdot < h_{2k_{0}}$ such that $h_{2j}=h_{2j-1}+2$ for $1\leq j \leq k_{0}$,
	$$\{h_{1}, h_{2}, \cdot\cdot\cdot, h_{2k_{0}}\}$$
is admissible and 
$$h_{2k_{0}}-h_{1} = O(k_{0}(\log k_{0})^{2}).$$
\end{lemma}

Before presenting the next lemma, we revisit the definition of ‘‘exponent of distribution to smooth moduli".

\begin{definition}\label{expsm}
	An arithmetic function $f$ with support contained in $[x, 2x)$ has exponent of distribution $\theta$ to smooth moduli if for every $\epsilon >0$ there exists a $\delta > 0$ for which the following holds.

	For any $P \in \{d \in \mathbb{N^{+}}: \mu(d) \neq 0, \ p | d \Rightarrow p \leq x^{\delta}\}$, any integer $a$ with $(a, P)=1$ and any $A > 0$ we have
	\[ \sum_{\substack{q \leq x^{\theta - \epsilon} \\ q | P}} \left| \sum_{n\equiv a\, (\!\bmod q)} f(n) - \frac{1}{\phi(q)}\sum_{(n,q)=1}f(n)\right| \ll_{\epsilon, A} x (\log x)^{-A}.
	\]
\end{definition}

We now present a minorant for the indicator function of primes, as constructed by Baker and Irving \cite{BI2017}. This minorant is equipped with a more robust equidistribution theorem in arithmetic progressions with smooth moduli.

\begin{lemma} \label{BIM} For all large $x$ there exists an arithmetic function $\rho(n)$ with support contained in $[x,2x)$ satisfying the following properties:
\begin{itemize}
\item [1.] $\rho(n)$ is a minorant for the indicator function of the primes, that is 
\begin{equation*}   \rho(n) \leq \left\{
	\begin{array}{rcl}
		1     &      & \mbox{$n$ is a prime}  \\
		0               &      & \mbox{otherwise.}
	\end{array} \right.
\end{equation*}
\item [2.] If $\rho(n) \neq 0$ then all prime factors of $n$ exceed $x^{\xi}$, for some fixed $\xi > 0$.
\item [3.] The function $\rho(n)$ has exponent of distribution $\theta$ to smooth moduli, where $\theta= \frac{1}{2}+\frac{7}{300}+\frac{17\eta}{120}$ for some $\eta \in (0,  \frac{22}{3295})$.
\item [4.] We have
\[ \sum_{x \leq n < 2x} \rho(n) = (1-c_{1}+o(1)) \frac{x}{\log x}
\]
for some $c_{1} < 8\times 10^{-6}$ such that $(1-c_{1})\theta > 0.52427$.
\end{itemize} 
\begin{proof}
	See R. C. Baker and A. J. Irving \cite[Lemmas 1, 2 and Section 5]{BI2017}.
\end{proof}
\end{lemma}

\section{Proof of Theorem \ref{NewTh}} \label{PMTHm}
\subsection{Setup} \label{SU}

\

Suppose that $x$ is sufficiently large. Let $\rho, \theta, c_{1}$, and $\xi$ be as in Lemma \ref{BIM}. We can choose $\epsilon$ sufficiently small, such that
	\begin{equation} \label{relcte}
		(1-c_{1})(\theta-\epsilon) > 0.52427.
	\end{equation}
Since $\rho$ has exponent of distribution $\theta$ to smooth moduli, we can find a $\delta > 0$ for which the following holds.

For any $P$ which is a product of distinct primes smaller than $x^{\delta}$, any integer $a$ with $(a, P)=1$ and any $A > 0$ 
\begin{equation} \label{rhosm}
	\sum_{\substack{q \leq x^{\theta - \epsilon} \\ q | P}} \left| \sum_{n\equiv a\, (\!\bmod q)} \rho(n) - \frac{1}{\phi(q)}\sum_{(n,q)=1}\rho(n)\right| \ll_{\epsilon, A} x (\log x)^{-A},
\end{equation}
 Let $\theta_{0}=\theta-\epsilon, R=x^{\theta_{0}/2-1/(100000m)}$ and
	\begin{equation} \label{Chk0}
	 k_{0}=m^{2}e^{\frac{4m}{\theta_{0}(1-c_{1})}+8}.
	 \end{equation}
We clearly have from \eqref{relcte} and Property 3 in Lemma \ref{BIM},
 \begin{equation} \label{snres}
 	\frac{4}{\theta_{0}(1-c_{1})} < \frac{4}{0.52427} < 7.63 \quad\quad  \mbox{and} \quad \quad 
 	\frac{1}{2} < \theta_{0} <  \frac{1}{2}+\frac{7}{300}+\frac{17}{120}\cdot\frac{22}{3295} = \frac{691}{1318}.
 \end{equation}
 
 Suppose that $\{h_{1}, \cdot\cdot\cdot h_{2k_{0}}\}$ is an admissible set described in Lemma \ref{adsetl}. First we use the $W$-trick. Set $W=\prod_{p<D_{0}}p$ for some $D_{0}$, by the Chinese remainder theorem, we can then find an integer $v$, such that $v+h_{i}$ is coprime to $W$ for each $h_{i}$. We restrict $n$ to be in this fixed residue class $v$ modulo $W$. One can choose $D_{0}=\log \log \log x$, so that $W \sim (\log \log x)^{1+o(1)}$ by an application of the prime number theorem. For a positive number $C$, we denote by $S(x,C)$ the quantity 
 \begin{equation*} 
 	\sum_{\substack{x \leq n < 2x \\ n\equiv v \, (\! \bmod  W) \\ \mu(n+h_{2i})\neq 0, 1\leq i \leq k_{0}}} \left( \sum_{j=1}^{k_{0}} \mathbf{1}_{\mathbb{P}}(n+h_{2j-1})\left( 1- \frac{\tau(n+h_{2j})}{C}\right)-m\right)\left( \sum_{\substack{d_{i} | n+h_{i}, 1\leq i \leq 2k_{0}	 }} \lambda_{\underline{d}}\right)^{2},
 \end{equation*}
 where $\lambda_{\underline{d}}$ are real constants to be chosen later.
 
 We wish to show, with an appropriate choice of $C$, that $S(x,C)>0$ for all large $x$. If $S(x,C)>0$ for some $x$, then at least one term in the sum over $n$ must have a strictly positive contribution. Since the sieve weights $\left( \sum \lambda_{d_{1}, \cdot\cdot\cdot, d_{2k_{0}}}\right)^{2}$ are nonnegative, we see that if there is a positive contribution from $n \in [x, 2x)$, then there exist distinct $1\leq j_{1}, \cdot\cdot\cdot, j_{m+1} \leq k_{0}$ such that
 \[ \mathbf{1}_{\mathbb{P}}(n+h_{2j_{i}-1})\left( 1- \frac{\tau(n+h_{2j_{i}})}{C}\right) > 0,
 \]
 i.e., $n+h_{2j_{i}-1}$ is prime and $\tau(n+h_{2j_{i}})< C$. Since $\mu(n+h_{2j_{i}})\neq 0$, we get that
 \[ \Omega(n+h_{2j_{i}})=\Omega(n+h_{2j_{i}-1}+2) \leq \frac{\log C}{\log 2}.
 \]
 Since this holds for all large $x$, there must be infinitely many integers $n$ such that $m+1$ elements of $(n+h_{2j-1})_{j=1}^{k_{0}}$ are prime and $n+h_{2j-1}+2$ has at most $\frac{\log C}{\log 2}$ prime factors. Furthermore, due to \eqref{Chk0} and \eqref{snres}, we have $h_{2k_{0}}-h_{1}=O(k_{0}(\log k_{0})^{2})= O(e^{7.63m})$ from Lemma \ref{adsetl}. Hence, Theorem \ref{NewTh} follows by showing $\log C \leq 7.36m + 4\log m + 21\log 2$.
 
 We shall choose $\lambda_{\underline{d}}$ in terms of a fixed symmetric function $f \! : [0,\infty)^{2k_{0}}\rightarrow \mathbb{R},$ supported on the truncated simplex 
 $$\Delta_{2k_{0}}^{[\kappa]}(1):= \{(t_{1},\cdot\cdot\cdot,t_{2k_{0}})\in [0,\kappa]^{2k_{0}}:t_{1}+\cdot\cdot\cdot +t_{2k_{0}}\leq 1\},$$ as 
 \begin{align*}
 	\lambda_{\underline{d}}=\mu(d_{1})\cdot\cdot\cdot \mu(d_{2k_{0}})f\left(\frac{\log d_{1}}{\log R}, \cdot\cdot\cdot ,\frac{\log d_{2k_{0}}}{\log R}\right),
 \end{align*}
 where $\kappa=2\min \{ \xi, \delta\}/\theta_{0}$. Hence, 
 \begin{align}\label{reswei}
 	d_{i} \leq R^{\kappa} \leq x^{\min \{ \xi, \delta\}}, \qquad \mbox{for}\  1\leq i \leq 2k_{0}
 \end{align}
provided $\lambda_{\underline{d}} \neq 0$.
 We further rewrite $S(x, C)$ as 
 \begin{align}\label{rewS}
 	S(x, C)= S_{1} - C^{-1}S_{2}-mS_{3},
 \end{align}
where
\[  S_{1} =\sum_{\substack{x \leq n < 2x \\ n\equiv v \, (\! \bmod  W) \\ \mu(n+h_{2i})\neq 0, 1\leq i \leq k_{0}}} \sum_{j=1}^{k_{0}} \mathbf{1}_{\mathbb{P}}(n+h_{2j-1})\left( \sum_{\substack{d_{i} | n+h_{i}, 1\leq i \leq 2k_{0}	 }} \lambda_{\underline{d}}\right)^{2},
\]
\[ S_{2}=\sum_{\substack{x \leq n < 2x \\ n\equiv v \, (\! \bmod  W) \\ \mu(n+h_{2i})\neq 0, 1\leq i \leq k_{0}}} \sum_{j=1}^{k_{0}}\mathbf{1}_{\mathbb{P}}(n+h_{2j-1})\tau(n+h_{2j})\left( \sum_{\substack{d_{i} | n+h_{i}, 1\leq i \leq 2k_{0}	 }} \lambda_{\underline{d}}\right)^{2},
\]
\[ S_{3}=\sum_{\substack{x \leq n < 2x \\ n\equiv v \, (\! \bmod  W) \\ \mu(n+h_{2i})\neq 0, 1\leq i \leq k_{0}}} \left( \sum_{\substack{d_{i} | n+h_{i}, 1\leq i \leq 2k_{0}	 }} \lambda_{\underline{d}}\right)^{2}.
\]
Based on the above discussion, the remainder of this paper is devoted to choosing a suitable function $f$ and estimating $S_{1}$, $S_{2}$, as well as $S_{3}$.

\subsection{The choice of the function $f$} \label{TCFf}

\

We set
\begin{equation*}
	\delta_{1}=\frac{1}{4.5k_{0}\log k_{0}}.
\end{equation*}
Let $h_{1}(t_{1},\cdot\cdot\cdot,t_{2k_{0}}) \!: [0,\infty)^{2k_{0}} \rightarrow \mathbb{R}$ be a smooth function with $| h_{1}(t_{1},\cdot\cdot\cdot,t_{2k_{0}}) | \leq1$ such that
\begin{equation*}  h_{1}(t_{1},\cdot\cdot\cdot,t_{2k_{0}})=\left\{
	\begin{array}{rcl}
		1,    &      &\mbox{if} \ \ (t_{1},\cdot\cdot\cdot,t_{2k_{0}})\in \Delta_{2k_{0}}(1-\delta_{1}),       \\
		0,              &      & \mbox{if}\ \  (t_{1},\cdot\cdot\cdot,t_{2k_{0}})\notin \Delta_{2k_{0}}(1).
	\end{array} \right.
\end{equation*}
Furthermore, we may assume that
\begin{align*} 
	\left| \frac{\partial h_{1}}{\partial t_{i}}(t_{1},\cdot\cdot\cdot,t_{2k_{0}})\right| \leq \frac{1}{\delta_{1}}+1
\end{align*}
for each $(t_{1},\cdot\cdot\cdot,t_{2k_{0}})\in \Delta_{2k_{0}}(1)\setminus \Delta_{2k_{0}}(1-\delta_{1})$ and $1\leq i \leq 2k_{0}$.
Let 
\[A= \log(2k_{0})-2\log \log (2k_{0})
\]
and
\[ T=\frac{e^{A}-1}{A}.
\]
It is obvious that for large $k_{0}$
\[ A > 0.99 \log k_{0}.
\]
Let
\[ \delta_{2}=\frac{\delta_{1}T}{10}.
\]
We also have 
\begin{equation*}
	\delta_{2} \geq \frac{1}{23(\log k_{0})^{4}}
\end{equation*}
for large $k_{0}$. Let $h^{*}_{2}(t) \! :[0,\infty)\rightarrow \mathbb{R}$ be a smooth function with $| h^{*}_{2}(t) | \leq 1$ such that
\begin{equation*}  h^{*}_{2}(t)=\left\{
	\begin{array}{rcl}
		1,     &      &\mbox{if} \ \ 0 \leq t\leq T-\delta_{2},       \\
		0,               &      & \mbox{if} \ \  t>T.
	\end{array} \right.
\end{equation*}
We may also assume that
\begin{align*} 
	\left| \frac{\mathrm{d} h^{*}_{2}}{\mathrm{d} t}(t) \right| \leq \frac{1}{\delta_{2}}+1
\end{align*}
for each $T-\delta_{2}\leq t \leq T$. 
Finally, we define\footnote{In Li and Pan's paper \cite{L-P2015},  the smooth function $f(t_{1},\cdot\cdot\cdot,t_{2k_{0}})$ is originally defined on $\mathbb{R}^{2k_{0}}$ rather than on $[0,\infty)^{2k_{0}}$. The authors additionally imposed the condition that $f(t_{1},\cdot\cdot\cdot,t_{2k_{0}})$ should vanish if $t_{i} < 0$ for some $1\leq i \leq 2k_{0}$. However, these restrictions can be relaxed when employing the Fourier analytic method to derive the Maynard-Tao sieve. For example, refer to \cite[Lemma 30]{PM2014} or \cite[Lemma 3.4]{A2018}.} the function $f \! : [0,\infty)^{2k_{0}} \rightarrow \mathbb{R} $ by
\begin{align*} 
	f(\underline{t})=\int_{t_{1}}^{\infty}\cdot\cdot\cdot \int_{t_{2k_{0}}}^{\infty}h_{1}(\underline{t})\prod_{i=1}^{2k_{0}}\frac{h^{*}_{2}(2k_{0}t_{i})}{1+2k_{0}At_{i}}\:\mathrm{d} \underline{t},\ \ \ \mbox{for} \ \underline{t}\in [0,\infty)^{2k_{0}}.
\end{align*}
As $h_{1}(\underline{t})\prod_{j=1}^{2k_{0}}\frac{h_{2}(2k_{0}t_{j})}{1+2k_{0}At_{j}}$ is a smooth function supported on $\Delta_{2k_{0}}^{[\frac{T}{2k_{0}}]}(1)$, we obtain that $f(\underline{t})$ is also a smooth function supported on $\Delta_{2k_{0}}^{[\frac{T}{2k_{0}}]}(1)$ and 
\begin{equation} \label{deF} 
	\frac{\partial^{2k_{0}}f(t_{1},\cdot\cdot\cdot,t_{2k_{0}})}{\partial t_{1}\cdot\cdot\cdot \partial t_{2k_{0}}}=h_{1}(\underline{t})\prod_{i=1}^{2k_{0}}\frac{h^{*}_{2}(2k_{0}t_{i})}{1+2k_{0}At_{i}}.
\end{equation}
Note that $\frac{T}{2k_{0}}\sim (\log 2k_{0})^{-3} \leq \kappa$ for large $k_{0}$. Thus we have $\supp f \subseteq \Delta_{2k_{0}}^{[\kappa]}(1)$ provided $k_{0}$ is large.
\subsection{A lower bound for $S_{1}$} \label{ALBS1}

\

We first note that \cite[eq. (3.5)]{L-P2015} gives
 \begin{align} \label{LBPRS2}
 	&\sum_{\substack{x \leq n < 2x \\ n\equiv v \, (\! \bmod  W) \\ \mu(n+h_{i})\neq 0, 1\leq i \leq k_{0}}}  \mathbf{1}_{\mathbb{P}}(n+h_{2j-1})\left( \sum_{\substack{d_{i} | n+h_{i}\forall i}} \lambda_{\underline{d}}\right)^{2}  
 	\notag
 	\\
 	= &(1+o(1))
 	\sum_{\substack{x \leq n < 2x \\ n\equiv v \, (\! \bmod  W) }}  \mathbf{1}_{\mathbb{P}}(n+h_{2j-1})\left( \sum_{\substack{d_{i} | n+h_{i}\forall i}} \lambda_{\underline{d}}\right)^{2}.
 \end{align}
We replace $\mathbf{1}_{\mathbb{P}}$ by $\rho$, expand out the square, and swap the order of summation to give 
	 \begin{align} \label{reprou}
		 \sum_{\substack{x \leq n < 2x \\ n\equiv v \, (\! \bmod  W) }}  \mathbf{1}_{\mathbb{P}}(n+h_{2j-1})\left( \sum_{\substack{d_{i} | n+h_{i}\forall i}} \lambda_{\underline{d}}\right)^{2}  \geq \sum_{\substack{x \leq n < 2x \\ n\equiv v \, (\! \bmod  W) }}  \rho(n+h_{2j-1})\left( \sum_{\substack{d_{i} | n+h_{i} \forall i}} \lambda_{\underline{d}}\right)^{2}
		\notag
		\\
		=  \sum_{\substack{d_{1}, \cdot\cdot\cdot, d_{2k_{0}}  \\ e_{1}, \cdot\cdot\cdot, e_{2k_{0}}}} \lambda_{d_{1}, \cdot\cdot\cdot, d_{2k_{0}}} \lambda_{e_{1}, \cdot\cdot\cdot, e_{2k_{0}}} \sum_{\substack{x \leq n < 2x \\ n\equiv v \, (\! \bmod  W) \\ [d_{i}, e_{i}] | n+h_{i} \forall i }} \rho(n+h_{2j-1}).
	\end{align}
 By the Chinese remainder theorem, the inner sum can be written as a sum over a single residue class modulo $q=W\prod_{i=1}^{2k_{0}}[d_{i},e_{i}]$, provided that $W, [d_{1}, e_{1}], \cdot\cdot\cdot, [d_{2k_{0}}, e_{2k_{0}}]$ are pairwise coprime. The integer $n+h_{2j-1}$ will lie in a residue class coprime to the modulus if and only if $d_{2j-1}=e_{2j-1}=1$. In this case, the inner sum will contribute $\frac{1}{\phi(q)}\sum_{(n,q)=1}\rho(n)+E(x, q, a)$, where
 \begin{align*}
 	E(x,q, a)= & \left| \sum_{\substack{x+h_{2j-1} \leq n < 2x+h_{2j-1} \\ n\equiv a \, (\! \bmod  q) }} \rho(n) - \frac{1}{\phi(q)}\sum_{\substack{ x\leq n < 2x \\ (n,q)=1 }}\rho(n)\right|
 	\notag
 	\\
 	= & \left| \sum_{\substack{x \leq n < 2x \\ n\equiv a \, (\! \bmod  q) }} \rho(n) - \frac{1}{\phi(q)}\sum_{\substack{x\leq n < 2x \\ (n,q)=1}}\rho(n)\right|+O(1),
 \end{align*}
and $a$ may depend on $W, d_{1}, \cdot\cdot\cdot, d_{2j-2}, d_{2j},\cdot\cdot\cdot,d_{2k_{0}}, e_{1}, \cdot\cdot\cdot, e_{2j-2}, e_{2j},\cdot\cdot\cdot,e_{2k_{0}}$.
Since $v+h_{i}$ is coprime to $W$ for each $h_{i}$,  $|h_{i}-h_{i'}| < D_{0}$ for all distinct $i,j$, due to \eqref{reswei} and since $\rho$ satisfies Property 2 in Lemma \ref{BIM}, the contribution of the inner sum in \eqref{reprou} is zero if either one pair of  $W, [d_{1}, e_{1}], \cdot\cdot\cdot, [d_{2k_{0}}, e_{2k_{0}}]$ share a common factor, or if either $d_{2j-1}$ or $e_{2j-1}$ are not 1. Thus we obtain
 \begin{align} \label{poc}
	&\sum_{\substack{x \leq n < 2x \\ n\equiv v \, (\! \bmod  W) }}  \mathbf{1}_{\mathbb{P}}(n+h_{2j-1})\left( \sum_{\substack{d_{i} | n+h_{i}\forall i}} \lambda_{\underline{d}}\right)^{2}  
	\notag
	\\
	\geq  &\sideset{}{'}{\sum}_{\substack{d_{1}, \cdot\cdot\cdot, d_{2k_{0}} \\ e_{1}, \cdot\cdot\cdot, e_{2k_{0}} \\ d_{2j-1}=e_{2j-1}=1}} \frac{\lambda_{\underline{d} }\lambda_{\underline{e}}}{\phi(q)} \sum_{\substack{(n,q)=1 \\ x\leq n < 2x}}\rho(n)-O\left(\sideset{}{'}{\sum}_{\substack{d_{1}, \cdot\cdot\cdot, d_{2k_{0}} \\ e_{1}, \cdot\cdot\cdot, e_{2k_{0}} \\ d_{2j-1}=e_{2j-1}=1}} \left|\lambda_{\underline{d} }\lambda_{\underline{e}}E(x,q,a)\right|\right)
	\notag
	\\
	= & \sideset{}{'}{\sum}_{\substack{d_{1}, \cdot\cdot\cdot, d_{2k_{0}} \\ e_{1}, \cdot\cdot\cdot, e_{2k_{0}} \\ d_{2j-1}=e_{2j-1}=1}} \frac{\lambda_{\underline{d} }\lambda_{\underline{e}}}{\phi(q)} \sum_{x\leq n < 2x}\rho(n)-O\left(\sideset{}{'}{\sum}_{\substack{d_{1}, \cdot\cdot\cdot, d_{2k_{0}} \\ e_{1}, \cdot\cdot\cdot, e_{2k_{0}} \\ d_{2j-1}=e_{2j-1}=1}} \left|\lambda_{\underline{d} }\lambda_{\underline{e}}E(x,q,a)\right|\right),
\end{align}
where $\sideset{}{'}{\sum}$ is used to denote the restriction that we require $W, [d_{1}, e_{1}], \cdot\cdot\cdot, [d_{2k_{0}}, e_{2k_{0}}]$ to be pairwise coprime. We can remove the restriction $(n,q)=1$ in the last step in view of \eqref{reswei} and $\rho$ since satisfies Property 2 in Lemma \ref{BIM}. 
Set
 $$F(\underline{t}) =\frac{\partial^{2k_{0}}f(t_{1},\cdot\cdot\cdot,t_{2k_{0}})}{\partial t_{1}\cdot\cdot\cdot \partial t_{2k_{0}}}.$$
By invoking $\sum_{x\leq n < 2x}\rho(n)=(1-c_{1}+o(1))x/ \log x$ and applying \cite[Lemma 3.4]{A2018} to the first sum in \eqref{poc}, we obtain (cf. \cite[Lemma 4.3]{A2018}) a main term of
\[ \frac{(1-c_{1} +o(1))x}{ (\log R)^{2k_{0}-1}\log x} \cdot \frac{W^{2k_{0}-1}}{\phi(W)^{2k_{0}}} \int_{\Delta_{2k_{0}-1}(1)} \left( \int_{0}^{1} F(\underline{t}) \:\mathrm{d}t_{2j-1} \right)^{2}\:\mathrm{d}t_{1} \cdot\cdot\cdot \mathrm{d}t_{2j-2}\:\mathrm{d}t_{2j}\cdot\cdot\cdot \mathrm{d}t_{2k_{0}-1}.
\]

Following the same argument as in \cite[p. 23, Subsection: The Motohashi-Pintz-Zhang case]{PM2014} along with \eqref{rhosm} and \eqref{reswei}, one can show the error term in \eqref{poc} contributes $\ll x(\log x)^{-A}$ for any fixed $A>0$. Combining this with the symmetry of $F$, we deduce that
\begin{align} \label{LBP}
	&\sum_{\substack{x \leq n < 2x \\ n\equiv v \, (\! \bmod  W) }}  \mathbf{1}_{\mathbb{P}}(n+h_{2j-1})\left( \sum_{\substack{d_{i} | n+h_{i}\forall i}} \lambda_{\underline{d}}\right)^{2} 
	\notag
	\\
	 \geq
	&\frac{(1-c_{1} +o(1))x}{ (\log R)^{2k_{0}-1}\log x} \cdot \frac{W^{2k_{0}-1}}{\phi(W)^{2k_{0}}} \int_{\Delta_{2k_{0}-1}(1)} \left( \int_{0}^{1} F(\underline{t}) \:\mathrm{d}t_{2{k_{0}}} \right)^{2}\:\mathrm{d}t_{1}\cdot\cdot\cdot \mathrm{d}t_{2k_{0}-1}.
\end{align}

Let 
\begin{equation*}
	\gamma =\frac{1}{A}\left(1-\frac{1}{1+AT}\right).
\end{equation*}
According to Maynard's work (cf. \cite[eq. (7.4) and (7.21)]{J2015}), we have 
\begin{align} \label{MWU}
	\int_{\Delta_{2k_{0}}(1)}  F^{\circ}(t_{1}, \cdot\cdot\cdot, t_{2k_{0}}) ^{2}\:\mathrm{d}t_{1} \cdot\cdot\cdot  \mathrm{d}t_{2k_{0}} \leq \frac{\gamma^{2k_{0}}}{(2k_{0})^{2k_{0}}}
\end{align}
and
\begin{align}\label{MWL}
	\int_{\Delta_{2k_{0}-1}(1)} \left( \int_{0}^{1} F^{\circ}(\underline{t}) \:\mathrm{d}t_{2{k_{0}}} \right)^{2}\:\mathrm{d}t_{1}\cdot\cdot\cdot \mathrm{d}t_{2k_{0}-1} \geq \frac{\log(2k_{0})-2\log \log (2k_{0})-2}{2k_{0}}\cdot\frac{\gamma^{2k_{0}}}{(2k_{0})^{2k_{0}}},
\end{align}
where
\[
F^{\circ}(t_{1} \cdot\cdot\cdot t_{2k_{0}})=\mathbf{1}_{\Delta_{2k_{0}}(1)}(t_{1} \cdot\cdot\cdot t_{2k_{0}})\prod_{j=1}^{2k_{0}}\frac{\mathbf{1}_{[0,T]}(2k_{0}t_{j})}{1+2k_{0}At_{j}}.
\]
On the other hand, \cite [eq. (3.11)]{L-P2015} gives\footnote{We use a slightly different notation: the function $F^{*}$ in \cite [eq. (3.11)]{L-P2015} corresponds to $F$ in this context.}
\begin{align}\label{chF}
	&\int_{\Delta_{2k_{0}-1}(1)} \left( \int_{0}^{1} F(\underline{t}) \:\mathrm{d}t_{1} \right)^{2}\:\mathrm{d}t_{2}\cdot\cdot\cdot \mathrm{d}t_{2k_{0}-1} 
	\notag
	\\
	\geq &(1-2.24k_{0}\delta_{1}) \int_{\Delta_{2k_{0}-1}(1)} \left( \int_{0}^{1} F^{\circ}(\underline{t}) \:\mathrm{d}t_{2{k_{0}}} \right)^{2}\:\mathrm{d}t_{1}\cdot\cdot\cdot \mathrm{d}t_{2k_{0}-1}.
\end{align}
It is easy to verify that
\[ 1-2.24k_{0}\delta_{1}=1-\frac{2.24k_{0}}{4.5k_{0}\log k_{0}}\geq \frac{\log(2k_{0})-2\log \log (2k_{0})-2.5}{\log(2k_{0})-2\log \log (2k_{0})-2}.
\]
Combining this with \eqref{MWL} and \eqref{chF}, we obtain 
\begin{align} \label{pFF*}
	\int_{\Delta_{2k_{0}-1}(1)} \left( \int_{0}^{1} F(\underline{t}) \:\mathrm{d}t_{2{k_{0}}} \right)^{2}\:\mathrm{d}t_{1}\cdot\cdot\cdot \mathrm{d}t_{2k_{0}-1}
	\geq \frac{\log(2k_{0})-2\log \log (2k_{0})-2.5}{2k_{0}}\cdot\frac{\gamma^{2k_{0}}}{(2k_{0})^{2k_{0}}}.
\end{align}
Since $m$ is sufficiently large, we have from \eqref{Chk0} and \eqref{snres}
\[ \log (2k_{0}) -2 \log \log (2k_{0}) \geq \frac{4m}{\theta_{0}(1-c_{1})} +4.628.
\]
Combining this with \eqref{LBPRS2}, \eqref{LBP}, and \eqref{pFF*}, we arrive at the following lower bound for $S_{1}$:
\begin{align} \label{s1finaLBPRS}
	S_{1} \geq \frac{(k_{0}(1-c_{1})+o(1))x}{(\log R)^{2k_{0}-1} \log x} \frac{W^{2k_{0}-1}}{\phi(W)^{2k_{0}}}
	\left(\frac{4m}{\theta_{0}(1-c_{1})} +2.128\right)\frac{\gamma^{2k_{0}}}{(2k_{0})^{2k_{0}+1}}.
\end{align}

\subsection{An upper bound for $S_{3}$} \label{AUBS3}

\

According to \cite[Lemma 4.2]{A2018}, one has
\begin{align} \label{LBNP}
	\sum_{\substack{x \leq n < 2x \\ n\equiv v \, (\! \bmod  W) }} \left( \sum_{\substack{d_{i} | n+h_{i}\forall i}} \lambda_{\underline{d}}\right)^{2} 
	= 
	\frac{(1+o(1))x}{ (\log R)^{2k_{0}}} \cdot \frac{W^{2k_{0}-1}}{\phi(W)^{2k_{0}}} \int_{\Delta_{2k_{0}}(1)}  F(\underline{t})^{2} \:\mathrm{d}\underline{t}.
\end{align}
Notice that $F(\underline{t}) \leq F^{\circ}(\underline{t})$. We conclude that from \eqref{LBNP} and \eqref{MWU} 
\begin{align} \label{s3FinaLBNP}
	S_{3}\leq\sum_{\substack{x \leq n < 2x \\ n\equiv v \, (\! \bmod  W) }} \left( \sum_{\substack{d_{i} | n+h_{i}\forall i}} \lambda_{\underline{d}}\right)^{2} 
	\leq &
	\frac{(1+o(1))x}{ (\log R)^{2k_{0}}} \cdot \frac{W^{2k_{0}-1}}{\phi(W)^{2k_{0}}} \int_{\Delta_{2k_{0}}(1)}  F^{\circ}(\underline{t})^{2} \:\mathrm{d}t_{1} \cdot\cdot\cdot \:\mathrm{d}t_{2k_{0}}
	\notag
	\\
	\leq &
	\frac{(1+o(1))x}{ (\log R)^{2k_{0}}} \cdot \frac{W^{2k_{0}-1}}{\phi(W)^{2k_{0}}} \frac{\gamma^{2k_{0}}}{(2k_{0})^{2k_{0}}}.
\end{align}
\subsection{An upper bound for $S_{2}$} \label{AUBS2}

\

In this section, we will use Tao's approach to establish an upper bound for $S_{2}$. Recall 
\[S_{2}= \sum_{\substack{x \leq n < 2x \\ n\equiv v \, (\! \bmod  W) \\ \mu(n+h_{2i})\neq 0, 1\leq i \leq k_{0}}} \sum_{j=1}^{k_{0}} \mathbf{1}_{\mathbb{P}}(n+h_{2j-1})\tau(n+h_{2j})\left( \sum_{\substack{d_{i} | n+h_{i}, 1\leq i \leq 2k_{0}	 }} \lambda_{\underline{d}}\right)^{2}.\]
From now on we only consider $j=1$. Let $\tilde{f}: [0,\infty)^{2k_{0}}\rightarrow \mathbb{R}$ be a smooth function with support on $\Delta_{2k_{0}}(\frac{2}{3\theta_{0}})$ such that $\tilde{f}(0,t_{2},\cdot\cdot\cdot,t_{2k_{0}})=f(0,t_{2},\cdot\cdot\cdot,t_{2k_{0}})$. Correspondingly, we define
\[
	\tilde{\lambda}_{\underline{d}}:=\mu(d_{1})\cdot\cdot\cdot \mu(d_{2k_{0}})\tilde{f}\left(\frac{\log d_{1}}{\log R}, \cdot\cdot\cdot ,\frac{\log d_{2k_{0}}}{\log R}\right).
\]
We therefore have 
 \begin{align} \label{ATT}
	&\sum_{\substack{x \leq n < 2x \\ n\equiv v \, (\! \bmod  W) \\ \mu(n+h_{2i})\neq 0, 1\leq i \leq k_{0}}}  \mathbf{1}_{\mathbb{P}}(n+h_{1})\tau(n+h_{2})\left( \sum_{\substack{d_{i} | n+h_{i}, 1\leq i \leq 2k_{0}	 }} \lambda_{\underline{d}}\right)^{2}
	\notag
	\\
	=&\sum_{\substack{x \leq n < 2x \\ n\equiv v \, (\! \bmod  W) \\ \mu(n+h_{2i})\neq 0, 1\leq i \leq k_{0}}}  \mathbf{1}_{\mathbb{P}}(n+h_{1})\tau(n+h_{2})\left( \sum_{\substack{d_{i} | n+h_{i}, 1\leq i \leq 2k_{0}	 }} \tilde{\lambda}_{\underline{d}}\right)^{2}
	\notag
	\\
	\leq& \sum_{\substack{x \leq n < 2x \\ n\equiv v \, (\! \bmod  W) }} \tau(n+h_{2})\left( \sum_{\substack{d_{i} | n+h_{i}, 1\leq i \leq 2k_{0}	 }} \tilde{\lambda}_{\underline{d}}\right)^{2}
	\notag
	\\
	=&\frac{x}{(\log R)^{2k_{0}}}\frac{W^{2k_{0}-1}}{\phi(W)^{2k_{0}}}\left(\frac{\log x}{\log R} \alpha (\tilde{f})- \beta_{1}(\tilde{f})-4\beta_{2}(\tilde{f}) + o(1) \right),
\end{align}
where
$$\alpha (\tilde{f})=\int_{\Delta_{2k_{0}}(\frac{2}{3\theta_{0}})}t_{2} \left(\frac{\partial^{2k_{0}+1}\tilde{f}(t_{1},\cdot\cdot\cdot,t_{2k_{0}})}{\partial t_{1}(\partial t_{2})^{2}\cdot\cdot\cdot \partial t_{2k_{0}}} \right)^{2}\:\mathrm{d} \underline{t},$$
$$\beta_{1}(\tilde{f})=\int_{\Delta_{2k_{0}}(\frac{2}{3\theta_{0}})}t_{2}^{2} \left(\frac{\partial^{2k_{0}+1}\tilde{f}(t_{1},\cdot\cdot\cdot,t_{2k_{0}})}{\partial t_{1}(\partial t_{2})^{2}\cdot\cdot\cdot \partial t_{2k_{0}}} \right)^{2}\:\mathrm{d} \underline{t},$$
and
$$\beta_{2}(\tilde{f})=\int_{\Delta_{2k_{0}}(\frac{2}{3\theta_{0}})}t_{2}\frac{\partial^{2k_{0}+1}\tilde{f}(t_{1},\cdot\cdot\cdot,t_{2k_{0}})}{\partial t_{1}(\partial t_{2})^{2}\cdot\cdot\cdot \partial t_{2k_{0}}} \frac{\partial^{2k_{0}}\tilde{f}(t_{1},\cdot\cdot\cdot,t_{2k_{0}})}{\partial t_{1}\cdot\cdot\cdot \partial t_{2k_{0}}} \:\mathrm{d} \underline{t}.$$
In the last step, we applied \cite[Lemma 5.10]{M-A2017} with $\mathcal{F}(t_{1},\cdot\cdot\cdot, t_{2k_{0}})=\tilde{f}(\frac{2}{3\theta_{0}}t_{1},\cdot\cdot\cdot, \frac{2}{3\theta_{0}}t_{2k_{0}})$ and $R=x^{\frac{1}{3}-\frac{2}{3}\cdot\frac{1}{10^{5}m\theta_{0}}}$, and then proceeded to change the variables in the integral. We will now carefully select the function $\tilde{f}$ to minimize $\alpha (\tilde{f})$ as much as possible since the main contribution will come from this term. 

By a converse to Cauchy-Schwarz and the fundamental theorem of calculus,
\begin{align} \label{CCS}
	\alpha (\tilde{f})&=\int_{\Delta_{2k_{0}-1}(\frac{2}{3\theta_{0}})}t_{2}\:\mathrm{d}t_{2}\cdot\cdot\cdot \mathrm{d}t_{2k_{0}} \int_{0}^{\frac{2}{3\theta_{0}}-t_{2}-\cdot\cdot\cdot-t_{2k_{0}}} \left(\frac{\partial^{2k_{0}+1}\tilde{f}(t_{1},\cdot\cdot\cdot,t_{2k_{0}})}{\partial t_{1}(\partial t_{2})^{2}\cdot\cdot\cdot \partial t_{2k_{0}}} \right)^{2}\:\mathrm{d} t_{1}
	\notag
	\\
	& \geq \int_{\Delta_{2k_{0}-1}(\frac{2}{3\theta_{0}})}\frac{t_{2}\:\mathrm{d}t_{2}\cdot\cdot\cdot \mathrm{d}t_{2k_{0}}}{\frac{2}{3\theta_{0}}-t_{2}-\cdot\cdot\cdot-t_{2k_{0}}}
	 \left(\int_{0}^{\frac{2}{3\theta_{0}}-t_{2}-\cdot\cdot\cdot-t_{2k_{0}}} \frac{\partial^{2k_{0}+1}\tilde{f}(t_{1},\cdot\cdot\cdot,t_{2k_{0}})}{\partial t_{1}(\partial t_{2})^{2}\cdot\cdot\cdot \partial t_{2k_{0}}} \:\mathrm{d} t_{1}\right)^{2}
	\notag
	\\
	&= \int_{\Delta_{2k_{0}-1}(\frac{2}{3\theta_{0}})}\frac{t_{2}}{\frac{2}{3\theta_{0}}-t_{2}-\cdot\cdot\cdot-t_{2k_{0}}}\left(\frac{\partial^{2k_{0}}f(0,t_{2},\cdot\cdot\cdot,t_{2k_{0}})}{(\partial t_{2})^{2}\cdot\cdot\cdot \partial t_{2k_{0}}} \right)^{2}\:\mathrm{d}t_{2}\cdot\cdot\cdot \mathrm{d}t_{2k_{0}}.
\end{align}
Moreover, the equality in \eqref{CCS} holds if and only if 
\[\frac{\partial^{2k_{0}+1}\tilde{f}(t_{1},\cdot\cdot\cdot,t_{2k_{0}})}{\partial t_{1}(\partial t_{2})^{2}\cdot\cdot\cdot \partial t_{2k_{0}}}=\frac{-\partial^{2k_{0}}f(0, t_{2}, \cdot\cdot\cdot,t_{2k_{0}})}{(\partial t_{2})^{2}\cdot\cdot\cdot \partial t_{2k_{0}}}\frac{1}{\frac{2}{3\theta_{0}}-t_{2}-\cdot\cdot\cdot-t_{2k_{0}}}.
\]
Equivalently,
\begin{equation} \label{idealc}
\frac{\partial^{2k_{0}}\tilde{f}(t_{1},\cdot\cdot\cdot,t_{2k_{0}})}{(\partial t_{2})^{2}\cdot\cdot\cdot \partial t_{2k_{0}}}=\frac{\partial^{2k_{0}}f(0, t_{2}\cdot\cdot\cdot,t_{2k_{0}})}{(\partial t_{2})^{2}\cdot\cdot\cdot \partial t_{2k_{0}}}\frac{\frac{2}{3\theta_{0}}-t_{1}-t_{2}-\cdot\cdot\cdot-t_{2k_{0}}}{\frac{2}{3\theta_{0}}-t_{2}-\cdot\cdot\cdot-t_{2k_{0}}}.
\end{equation}
Notice that the function on the right side of \eqref{idealc} is not supported on $\Delta_{2k_{0}}(\frac{2}{3\theta_{0}})$. Therefore, it is impossible for  $\alpha(\tilde{f})$ to attain the value on the right-hand side of \eqref{CCS}. However, we can still choose an appropriate function $\tilde{f}$ in such a way that $\alpha(\tilde{f})$  does not deviate too much from this ideal value.

Specifically, we define 
\begin{equation*} L(t_{1},\cdot\cdot\cdot,t_{2k_{0}}):=\left\{
	\begin{array}{rcl}
		\frac{\partial^{2k_{0}}f(0, t_{2}\cdot\cdot\cdot,t_{2k_{0}})}{(\partial t_{2})^{2}\cdot\cdot\cdot \partial t_{2k_{0}}}\frac{\frac{2}{3\theta_{0}}-t_{1}-t_{2}-\cdot\cdot\cdot-t_{2k_{0}}}{\frac{2}{3\theta_{0}}-t_{2}-\cdot\cdot\cdot-t_{2k_{0}}},    &      &\mbox{if} \ \ t_{2}+\cdot\cdot\cdot+t_{2k_{0}}\neq \frac{2}{3\theta_{0}},       \\
		0,              &      & \mbox{if}\ \  t_{2}+\cdot\cdot\cdot+t_{2k_{0}}= \frac{2}{3\theta_{0}}.
	\end{array} \right.
\end{equation*}
Let $h(t_{1},\cdot\cdot\cdot,t_{2k_{0}}) \!: [0,\infty)^{2k_{0}} \rightarrow \mathbb{R}$ be a smooth function with $| h(t_{1},\cdot\cdot\cdot,t_{2k_{0}}) | \leq1$ such that
\begin{equation} \label{smfh}  h(t_{1},\cdot\cdot\cdot,t_{2k_{0}})=\left\{
	\begin{array}{rcl}
		1,    &      &\mbox{if} \ \ (t_{1},\cdot\cdot\cdot,t_{2k_{0}})\in \Delta_{2k_{0}}(\frac{2}{3\theta_{0}}-\delta'),       \\
		0,              &      & \mbox{if}\ \  (t_{1},\cdot\cdot\cdot,t_{2k_{0}})\notin \Delta_{2k_{0}}(\frac{2}{3\theta_{0}}),
	\end{array} \right.
\end{equation}
where $\delta'$ is a small constant to be chosen soon. Furthermore, we can assume that
\begin{align} \label{deh} 
	\left| \frac{\partial h_{1}}{\partial t_{i}}(t_{1},\cdot\cdot\cdot,t_{2k_{0}})\right| \leq \frac{1}{\delta'}+1
\end{align}
for each $(t_{1},\cdot\cdot\cdot,t_{2k_{0}})\in \Delta_{2k_{0}}(\frac{2}{3\theta_{0}})\setminus \Delta_{2k_{0}}(\frac{2}{3\theta_{0}}-\delta')$ and $1\leq i \leq 2k_{0}$.
Finally, we select the function $\tilde{f}$ by
\[
\tilde{f}(\underline{t})=(-1)^{2k_{0}} \int_{t_{2}}^{\infty}\cdot\cdot\cdot \int_{t_{2k_{0}}}^{\infty}\left(\int_{y}^{\infty}h(t_{1},u_{2},\cdot\cdot\cdot,u_{2k_{0}})L(t_{1},u_{2},\cdot\cdot\cdot,u_{2k_{0}})\:\mathrm{d}u_{2}\right)\:\mathrm{d}y\mathrm{d}u_{3}\cdot\cdot\cdot \mathrm{d}u_{2k_{0}}.
\]
We clearly have 
$$\supp \tilde{f}\subseteq \Delta_{2k_{0}}(\frac{2}{3\theta_{0}}).$$
Note that $\supp f(0, t_{2},\cdot\cdot\cdot,t_{2k_{0}})\subseteq \Delta_{2k_{0}-1}(1)$. We also have $\tilde{f}(0,t_{2},\cdot\cdot\cdot,t_{2k_{0}})=f(0,t_{2},\cdot\cdot\cdot,t_{2k_{0}})$. Moreover, it follows from \eqref{smfh} and \eqref{deh} that
\begin{align} \label{diffnewf}
	\alpha(\tilde{f})&=\int_{\Delta_{2k_{0}}(\frac{2}{3\theta_{0}})}t_{2}\left(\frac{\partial}{\partial t_{1}}h(t_{1},\cdot\cdot\cdot,t_{2k_{0}})L(t_{1},\cdot\cdot\cdot,t_{2k_{0}})\right)^{2}\:\mathrm{d}\underline{t}
	\notag
	\\
	&=\int_{\Delta_{2k_{0}}(\frac{2}{3\theta_{0}}-\delta')}t_{2}\left(\frac{\partial L}{\partial t_{1}}\right)^{2}\:\mathrm{d}\underline{t}+\int_{\Delta_{2k_{0}}(\frac{2}{3\theta_{0}})\setminus \Delta_{2k_{0}}(\frac{2}{3\theta_{0}}-\delta')}t_{2}\left( \frac{\partial h}{\partial t_{1}} \cdot L + h \cdot \frac{\partial L}{\partial t_{1}} \right)^{2}\:\mathrm{d}\underline{t}
	\notag
	\\	
	& \leq \int_{\Delta_{2k_{0}}(\frac{2}{3\theta_{0}})}t_{2}\left(\frac{\partial L}{\partial t_{1}}\right)^{2}\:\mathrm{d}\underline{t}+2\int_{\Delta_{2k_{0}}(\frac{2}{3\theta_{0}})\setminus \Delta_{2k_{0}}(\frac{2}{3\theta_{0}}-\delta')}t_{2}\left( \frac{\partial h}{\partial t_{1}} \cdot L  \right)^{2}+t_{2}\left(h \cdot \frac{\partial L}{\partial t_{1}} \right)^{2}\:\mathrm{d}\underline{t}
	\notag
	\\
	&= \int_{\Delta_{2k_{0}}(\frac{2}{3\theta_{0}})}t_{2}\left(\frac{\partial L}{\partial t_{1}}\right)^{2}\:\mathrm{d}\underline{t} + O_{f}(\delta').
\end{align}
Here we used $L(\underline{t})\ll_{f} \delta'$ when $\underline{t} \in \Delta_{2k_{0}}(\frac{2}{3\theta_{0}})\setminus \Delta_{2k_{0}}(\frac{2}{3\theta_{0}}-\delta')$ and the volume of $ \Delta_{2k_{0}}(\frac{2}{3\theta_{0}})\setminus \Delta_{2k_{0}}(\frac{2}{3\theta_{0}}-\delta')$ is smaller than $\delta'$.

We now focus on the integral on the right-hand side of \eqref{diffnewf}, i.e.,
\begin{align} \label{mainalp}
	\int_{\Delta_{2k_{0}}(\frac{2}{3\theta_{0}})}t_{2}\left(\frac{\partial L}{\partial t_{1}}\right)^{2}\:\mathrm{d}\underline{t}= \int_{\Delta_{2k_{0}}(\frac{2}{3\theta_{0}})}t_{2}\left(\frac{-\partial^{2k_{0}}f(0, t_{2}, \cdot\cdot\cdot,t_{2k_{0}})}{(\partial t_{2})^{2}\cdot\cdot\cdot \partial t_{2k_{0}}}\frac{1}{\frac{2}{3\theta_{0}}-t_{2}-\cdot\cdot\cdot-t_{2k_{0}}}\right)^{2}\:\mathrm{d}\underline{t}
	\notag
	\\
	=\int_{\Delta_{2k_{0}-1}(\frac{2}{3\theta_{0}})}\frac{t_{2}}{\frac{2}{3\theta_{0}}-t_{2}-\cdot\cdot\cdot-t_{2k_{0}}}\left(\frac{\partial^{2k_{0}}f(0,t_{2},\cdot\cdot\cdot,t_{2k_{0}})}{(\partial t_{2})^{2}\cdot\cdot\cdot \partial t_{2k_{0}}} \right)^{2}\:\mathrm{d}t_{2}\cdot\cdot\cdot \mathrm{d}t_{2k_{0}} \quad \quad \quad \quad \quad 
	\notag
	\\
	=\int_{\Delta_{2k_{0}-1}(\frac{2}{3\theta_{0}})}\frac{t_{2}}{\frac{2}{3\theta_{0}}-t_{2}-\cdot\cdot\cdot-t_{2k_{0}}}\left(\int_{0}^{\infty}\frac{\partial^{2k_{0}+1}f(t_{1},t_{2},\cdot\cdot\cdot,t_{2k_{0}})}{\partial t_{1} (\partial t_{2})^{2}\cdot\cdot\cdot \partial t_{2k_{0}}}\:\mathrm{d}t_{1} \right)^{2}\:\mathrm{d}t_{2}\cdot\cdot\cdot \mathrm{d}t_{2k_{0}}.
	\end{align}
Recall that (cf. \eqref{deF})
\[ \frac{\partial^{2k_{0}}f(t_{1},\cdot\cdot\cdot,t_{2k_{0}})}{\partial t_{1}\cdot\cdot\cdot \partial t_{2k_{0}}}=h_{1}(\underline{t})\prod_{i=1}^{2k_{0}}\frac{h^{*}_{2}(2k_{0}t_{i})}{1+2k_{0}At_{i}}.
\]
We therefore have
\begin{align}\label{sepder}
\frac{\partial^{2k_{0}+1}f(t_{1},\cdot\cdot\cdot,t_{2k_{0}})}{\partial t_{1}(\partial t_{2})^{2}\cdot\cdot\cdot \partial t_{2k_{0}}}
=&\frac{\partial h_{1}}{\partial t_{2}}(\underline{t})\prod_{i=1}^{2k_{0}}\frac{h^{*}_{2}(2k_{0}t_{i})}{1+2k_{0}At_{i}}+h_{1}(\underline{t})\frac{2k_{0}h^{*'}_{2}(2k_{0}t_{2})}{1+2k_{0}At_{2}}\prod_{i \neq 2}\frac{h^{*}_{2}(2k_{0}t_{i})}{1+2k_{0}At_{i}}
\notag
\\
&-h_{1}(\underline{t_{1}})\frac{2k_{0}Ah^{*}_{2}(2k_{0}t_{2})}{(1+2k_{0}At_{2})^{2}}\prod_{i \neq 2}\frac{h^{*}_{2}(2k_{0}t_{i})}{1+2k_{0}At_{i}}.
\end{align}
Substituting  \eqref{sepder} into \eqref{mainalp} and applying Cauchy-Schwarz gives 
	\begin{align} \label{sub-cs-al}
		\int_{\Delta_{2k_{0}}(\frac{2}{3\theta_{0}})}t_{2}\left(\frac{\partial L}{\partial t_{1}}\right)^{2}\:\mathrm{d}\underline{t} \leq I_{1}+I_{2}+I_{3}+2\sqrt{I_{1}I_{2}}+2\sqrt{I_{1}I_{3}}+2\sqrt{I_{2}I_{3}},
	\end{align}
where
\[I_{1}=\int_{\Delta_{2k_{0}-1}(\frac{2}{3\theta_{0}})}\frac{t_{2}}{\frac{2}{3\theta_{0}}-\sum_{i=2}^{2k_{0}}t_{i}}\left(\int_{0}^{\infty}\frac{\partial h_{1}}{\partial t_{2}}(\underline{t})\prod_{i=1}^{2k_{0}}\frac{h^{*}_{2}(2k_{0}t_{i})}{1+2k_{0}At_{i}}\:\mathrm{d}t_{1} \right)^{2}\:\mathrm{d}t_{2}\cdot\cdot\cdot \mathrm{d}t_{2k_{0}},
\]
\[
I_{2}=\int_{\Delta_{2k_{0}-1}(\frac{2}{3\theta_{0}})}\frac{t_{2}}{\frac{2}{3\theta_{0}}-\sum_{i=2}^{2k_{0}}t_{i}}\left(\int_{0}^{\infty}h_{1}(\underline{t})\frac{2k_{0}h^{*'}_{2}(2k_{0}t_{2})}{1+2k_{0}At_{2}}\prod_{i \neq 2}\frac{h^{*}_{2}(2k_{0}t_{i})}{1+2k_{0}At_{i}}\:\mathrm{d}t_{1} \right)^{2}\:\mathrm{d}t_{2}\cdot\cdot\cdot \mathrm{d}t_{2k_{0}},
\]
\[
I_{3}=\int_{\Delta_{2k_{0}-1}(\frac{2}{3\theta_{0}})}\frac{t_{2}}{\frac{2}{3\theta_{0}}-\sum_{i=2}^{2k_{0}}t_{i}}\left(\int_{0}^{\infty}h_{1}(\underline{t})\frac{2k_{0}Ah^{*}_{2}(2k_{0}t_{2})}{(1+2k_{0}At_{2})^{2}}\prod_{i \neq 2}\frac{h^{*}_{2}(2k_{0}t_{i})}{1+2k_{0}At_{i}}\:\mathrm{d}t_{1} \right)^{2}\:\mathrm{d}t_{2}\cdot\cdot\cdot \mathrm{d}t_{2k_{0}}.
\]

We first deal with $I_{3}$. Noting that $\supp h_{1} \subseteq \Delta_{2k_{0}}(1)$, $h_{1}\leq 1$, $h^{*}_{2}\leq 1$ and $\supp h^{*}_{2} \subseteq [0, T]$, we have

\begin{align}\label{esI3}
	I_{3}\leq&\int_{\Delta_{2k_{0}-1}(1)}\frac{4k^{2}_{0}A^{2}t_{2}}{(1+2k_{0}At_{2})^{4}\left(\frac{2}{3\theta_{0}}-\sum_{i=2}^{2k_{0}}t_{i}\right)}\left(\int_{0}^{T/(2k_{0})}\frac{1}{1+2k_{0}At_{1}} \:\mathrm{d}t_{1}\right)^{2}\:\mathrm{d}t_{2}
	\notag
	\\
	 &\quad\quad \quad\quad \quad\quad\quad\quad \quad\quad \quad\quad\quad\quad \quad\quad \quad\quad\quad\quad \ \ \cdot \prod_{i \neq 1, 2}\frac{h^{*}_{2}(2k_{0}t_{i})^{2}\:\mathrm{d}t_{i}}{(1+2k_{0}At_{i})^{2}}
	\notag
	\\
	 =&\int_{\Delta_{2k_{0}-1}(1)}\frac{4k^{2}_{0}A^{2}t_{2}}{(1+2k_{0}At_{2})^{4}}\frac{1}{(2k_{0})^{2}\left(\frac{2}{3\theta_{0}}-\sum_{i=2}^{2k_{0}}t_{i}\right)} \:\mathrm{d}t_{2}
	  \prod_{i \neq 1, 2}\frac{h^{*}_{2}(2k_{0}t_{i})^{2}\:\mathrm{d}t_{i}}{(1+2k_{0}At_{i})^{2}}
	  \notag
	  \\
	  \leq&\left(\max_{0\leq r \leq 1} \frac{1}{\frac{2}{3\theta_{0}}-r} \right) \int_{0}^{\infty}\frac{A^{2}t_{2}\:\mathrm{d}t_{2}}{(1+2k_{0}At_{2})^{4}}\left(\int_{0}^{T/(2k_{0})}\frac{\mathrm{d}t}{(1+2k_{0}At)^2}\right)^{2k_{0}-2}
	  \notag
	  \\
	  =&\frac{1}{6} \cdot\frac{ 1}{(2k_{0})^{2}\left(\frac{2}{3\theta_{0}}-1\right)} \cdot \frac{\gamma^{2k_{0}-2}}{(2k_{0})^{2k_{0}-2}}.
\end{align}
Next, we have
\begin{align*}
	I_{2}\leq&\int_{\Delta_{2k_{0}-1}(1)}\frac{t_{2}}{\frac{2}{3\theta_{0}}-\sum_{i=2}^{2k_{0}}t_{i}}\left(\int_{0}^{T/(2k_{0})}\frac{1}{1+2k_{0}At_{1}} \:\mathrm{d}t_{1}\right)^{2}\left(\frac{2k_{0}h^{*'}_{2}(2k_{0}t_{2})}{1+2k_{0}At_{2}}\right)^{2}\:\mathrm{d}t_{2}
	\\
	&\quad\quad \quad\quad \quad\quad\quad\quad \quad\quad \quad\quad\quad\quad \quad\quad \quad\quad\quad\quad \quad \ \ \cdot \prod_{i \neq 1, 2}\frac{h^{*}_{2}(2k_{0}t_{i})^{2}\:\mathrm{d}t_{i}}{(1+2k_{0}At_{i})^{2}}
	\\
	=&\int_{\Delta_{2k_{0}-1}(1)}\frac{4k^{2}_{0}t_{2}h^{*'}_{2}(2k_{0}t_{2})^{2}}{(1+2k_{0}At_{2})^{2}}\frac{1}{(2k_{0})^{2}\left(\frac{2}{3\theta_{0}}-\sum_{i=2}^{2k_{0}}t_{i}\right)} \:\mathrm{d}t_{2}
	\prod_{i \neq 1, 2}\frac{h^{*}_{2}(2k_{0}t_{i})^{2}\:\mathrm{d}t_{i}}{(1+2k_{0}At_{i})^{2}}
	\\
	\leq& \left(\max_{0\leq r \leq 1} \frac{1}{\frac{2}{3\theta_{0}}-r} \right)\left(1+\frac{1}{\delta_{2}}\right)^{2} \int_{(T-\delta_{2})/(2k_{0})}^{T/(2k_{0})}\frac{t_{2}\:\mathrm{d}t_{2}}{(1+2k_{0}At_{2})^{2}}\left(\frac{\gamma}{2k_{0}}\right)^{2k_{0}-2}
	\\
	\leq& \frac{1}{\frac{2}{3\theta_{0}}-1} \left(1+\frac{1}{\delta_{2}}\right)^{2} \frac{\delta_{2}}{2k_{0}}\cdot\frac{T}{2k_{0}(1+(T-\delta_{2})A)^{2}}\cdot\frac{\gamma^{2k_{0}-2}}{(2k_{0})^{2k_{0}-2}}
	\\
	\leq& \frac{1}{\frac{2}{3\theta_{0}}-1}\left(1+\frac{1}{\delta_{2}}\right)^{2} \frac{\delta_{2}}{2k_{0}}\cdot\frac{T}{2k_{0}(1+AT)AT}\cdot\frac{\gamma^{2k_{0}-2}}{(2k_{0})^{2k_{0}-2}}
	\\
	\leq&\frac{1}{\frac{2}{3\theta_{0}}-1}\left(1+\frac{1}{\delta_{2}}\right)^{2} \frac{\delta_{2}}{2k_{0}}\cdot\frac{1.02 \log k_{0}}{(2k_{0})^{2}}\cdot\frac{\gamma^{2k_{0}-2}}{(2k_{0})^{2k_{0}-2}}
\end{align*}
by recalling that $A > 0.99 \log k_{0}$ and $1+AT=2k_{0}/§(\log (2k_{0}))^{2}$. Noting that $\delta_{2} \rightarrow 0$ as $k_{0} \rightarrow \infty$ and $\delta_{2} \geq \frac{1}{23(\log k_{0})^{4}}$, we arrive at as $k_{0} \rightarrow \infty$
\begin{align}\label{esI2}
	I_{2} = o\left(\frac{\gamma^{2k_{0}-2}}{(2k_{0})^{2k_{0}}}\right).
\end{align}
Finally, we turn to estimate $I_{1}$, i.e.,
\begin{align}\label{esI1}
	&\quad \ I_{1}=\int_{\Delta_{2k_{0}-1}(1)}\frac{t_{2}}{\frac{2}{3\theta_{0}}-\sum_{i=2}^{2k_{0}}t_{i}}\left(\int_{0}^{1-\sum_{i=2}^{2k_{0}}t_{i}}\frac{\partial h_{1}}{\partial t_{2}}(\underline{t})\frac{h^{*}_{2}(2k_{0}t_{1})\:\mathrm{d}t_{1}}{1+2k_{0}At_{1}}\right)^{2}
	\prod_{i \neq 1}\frac{h^{*}_{2}(2k_{0}t_{i})^{2}\:\mathrm{d}t_{i}}{(1+2k_{0}At_{i})^{2}} 
	\notag
	\\
	&=\int_{\Delta_{2k_{0}-1}(1-\delta_{1})}\frac{t_{2}}{\frac{2}{3\theta_{0}}-\sum_{i=2}^{2k_{0}}t_{i}}\left(\int_{1-\delta_{1}-\sum_{i=2}^{2k_{0}}t_{i}}^{1-\sum_{i=2}^{2k_{0}}t_{i}}\frac{\partial h_{1}}{\partial t_{2}}(\underline{t})\frac{h^{*}_{2}(2k_{0}t_{1})\:\mathrm{d}t_{1}}{1+2k_{0}At_{1}}\right)^{2}
	\prod_{i \neq 1}\frac{h^{*}_{2}(2k_{0}t_{i})^{2}\:\mathrm{d}t_{i}}{(1+2k_{0}At_{i})^{2}}+
	\notag
	\\
	&\int_{\Delta_{2k_{0}-1}(1)\setminus\Delta_{2k_{0}-1}(1-\delta_{1})}\frac{t_{2}}{\frac{2}{3\theta_{0}}-\sum_{i=2}^{2k_{0}}t_{i}}\left(\int_{0}^{1-\sum_{i=2}^{2k_{0}}t_{i}}\frac{\partial h_{1}}{\partial t_{2}}(\underline{t})\frac{h^{*}_{2}(2k_{0}t_{1})\:\mathrm{d}t_{1}}{1+2k_{0}At_{1}}\right)^{2}
    \prod_{i \neq 1}\frac{h^{*}_{2}(2k_{0}t_{i})^{2}\:\mathrm{d}t_{i}}{(1+2k_{0}At_{i})^{2}}
	\notag
	\\
	&\leq\int_{\Delta_{2k_{0}-1}(1-\delta_{1})}\frac{t_{2}}{\frac{2}{3\theta_{0}}-\sum_{i=2}^{2k_{0}}t_{i}}\left(1+\frac{1}{\delta_{1}}\right)^{2}\left(\frac{1}{2k_{0}A}\right)^{2}\left(\log \frac{1+2k_{0}A(1-\sum_{i=2}^{2k_{0}}t_{i})}{1+2k_{0}A(1-\delta_{1}-\sum_{i=2}^{2k_{0}}t_{i})}\right)^{2}
	\notag
	\\
	&\quad \quad \quad \quad \ \prod_{i \neq 1}\frac{h^{*}_{2}(2k_{0}t_{i})^{2}\:\mathrm{d}t_{i}}{(1+2k_{0}At_{i})^{2}}+\int_{\Delta_{2k_{0}-1}(1)\setminus\Delta_{2k_{0}-1}(1-\delta_{1})}\frac{t_{2}}{\frac{2}{3\theta_{0}}-\sum_{i=2}^{2k_{0}}t_{i}}\left(1+\frac{1}{\delta_{1}}\right)^{2}\left(\frac{1}{2k_{0}A}\right)^{2}
	\notag
	\\
	&\quad \quad \quad \quad \quad \quad\quad \quad\cdot \left(\log \left(1+2k_{0}A\left(1-\sum_{i=2}^{2k_{0}}t_{i}\right)\right)\right)^{2}\prod_{i \neq 1}\frac{h^{*}_{2}(2k_{0}t_{i})^{2}\:\mathrm{d}t_{i}}{(1+2k_{0}At_{i})^{2}}=: I_{1,1}+I_{1,2}.
\end{align}
Using the fact that
\[\log \frac{1+2k_{0}A(1-\sum_{i=2}^{2k_{0}}t_{i})}{1+2k_{0}A(1-\delta_{1}-\sum_{i=2}^{2k_{0}}t_{i})}\leq 2k_{0}A\delta_{1}
\]
for $\sum_{i=2}^{2k_{0}}t_{i} \leq 1-\delta_{1}$, we have
\begin{align}\label{esI11}
	I_{1,1}
	\leq&\int_{\Delta_{2k_{0}-1}(1-\delta_{1})}\frac{t_{2}}{\frac{2}{3\theta_{0}}-1}\left(1+\delta_{1}\right)^{2}
	\prod_{i \neq 1}\frac{h^{*}_{2}(2k_{0}t_{i})^{2}\:\mathrm{d}t_{i}}{(1+2k_{0}At_{i})^{2}}
	\notag
	\\
	\leq&\left(\frac{2}{3\theta_{0}}-1\right)^{-1}\left(1+\delta_{1}\right)^{2}\int_{0}^{T/2k_{0}}
	\frac{t_{2}\:\mathrm{d}t_{2}}{(1+2k_{0}At_{2})^{2}}\left(\int_{0}^{T/2k_{0}}
	\frac{\mathrm{d}t}{(1+2k_{0}At_{2})^{2}}\right)^{2k_{0}-2}
	\notag
	\\
	=&\left(\frac{2}{3\theta_{0}}-1\right)^{-1}\left(1+\delta_{1}\right)^{2}\frac{1-\gamma}{(2k_{0})^{2}A}\left(\frac{\gamma}{2k_{0}}\right)^{2k_{0}-2}.
\end{align}
Noting that 
\[\frac{t}{(1+2k_{0}At)^{2}}\leq \frac{1}{8k_{0}A}\]
for $t\geq 0$ and letting $r=t_{2}+\cdot\cdot\cdot+t_{2k_{0}}$, we have
\begin{align}\label{esI12}
	I_{1,2}
	\leq&  \frac{\left( \log (1+ 2k_{0}A\delta_{1})\right)^{2}}{(2k_{0}A)^{2}\left(\frac{2}{3\theta_{0}}-1\right)}\left(1+\frac{1}{\delta_{1}}\right)^{2} \int_{\Delta_{2k_{0}-1}(1)\setminus\Delta_{2k_{0}-1}(1-\delta_{1})}t_{2}
	\prod_{i \neq 1}\frac{h^{*}_{2}(2k_{0}t_{i})^{2}\:\mathrm{d}t_{i}}{(1+2k_{0}At_{i})^{2}}
	\notag
	\\
	\leq&\frac{\left( \log (1+ 2k_{0}A\delta_{1})\right)^{2}}{(2k_{0}A)^{2}\left(\frac{2}{3\theta_{0}}-1\right)}\left(1+\frac{1}{\delta_{1}}\right)^{2} \int_{\Delta_{2k_{0}-2}(1)}\int_{1 - \delta_{1}}^{1}\frac{|r-\sum_{i=3}^{2k_{0}}t_{i}|\:\mathrm{d}r}{(1+2k_{0}A|r-\sum_{i=3}^{2k_{0}}t_{i}|)^{2}}
	\notag
	\\
	&\quad\quad\quad\quad\quad\quad\quad\quad\quad\quad\quad\quad\quad\quad\quad\quad\quad\quad\quad\quad \ \cdot\prod_{i \neq 1, 2}\frac{h^{*}_{2}(2k_{0}t_{i})^{2}\:\mathrm{d}t_{i}}{(1+2k_{0}At_{i})^{2}}
	\notag
	\\
	\leq&\frac{\left( \log (1+ 2k_{0}A\delta_{1})\right)^{2}}{(2k_{0}A)^{2}\left(\frac{2}{3\theta_{0}}-1\right)}\left(1+\frac{1}{\delta_{1}}\right)^{2}\frac{\delta_{1}}{8k_{0}A}\left(\frac{\gamma}{2k_{0}}\right)^{2k_{0}-2}.
\end{align}
Note that $\gamma\rightarrow 0$  as $k_{0}\rightarrow \infty$ and $\delta_{1}= \frac{1}{4.5k_{0}\log k_{0}}$. From \eqref{esI1},  \eqref{esI11} and \eqref{esI12} we conclude that as $k_{0} \rightarrow \infty$
\begin{align}\label{esl1f}
	I_{1}\leq\frac{\left(1+\delta_{1}\right)^{2}(1-\gamma)}{\left(\frac{2}{3\theta_{0}}-1\right)(2k_{0})^{2}A}\left(\frac{\gamma}{2k_{0}}\right)^{2k_{0}-2}+\frac{\left( \log (1+ 2k_{0}A\delta_{1})\right)^{2}}{(2k_{0}A)^{2}\left(\frac{2}{3\theta_{0}}-1\right)}\left(1+\frac{1}{\delta_{1}}\right)^{2}\frac{\delta_{1}}{8k_{0}A}\left(\frac{\gamma}{2k_{0}}\right)^{2k_{0}-2}
	\notag
	\\
	=o\left(\frac{\gamma^{2k_{0}-2}}{(2k_{0})^{2k_{0}}}\right).\quad\quad\quad\quad\quad\quad\quad\quad\quad\quad\quad\quad\quad\quad\quad\quad\quad\quad\quad\quad\quad\quad\quad\quad\quad\quad\ \ 
\end{align}	
A combination of \eqref{diffnewf}, \eqref{sub-cs-al}, \eqref{esI3}, \eqref{esI2} and \eqref{esl1f} leads to
\begin{align*}
	\alpha(\tilde{f})&\leq O_{f}(\delta')+\left(\frac{1}{6} +o(1) \right)\frac{1}{\left(\frac{2}{3\theta_{0}}-1\right)}\frac{\gamma^{2k_{0}-2}}{(2k_{0})^{2k_{0}}}.
\end{align*}
Combining this with the choice of a sufficiently small $\delta'$ (depending on $k_{0}$), we conclude that for large $k_{0}$,
\begin{align}\label{esarp}
	\alpha(\tilde{f})\leq\frac{0.167}{\left(\frac{2}{3\theta_{0}}-1\right)}\frac{\gamma^{2k_{0}-2}}{(2k_{0})^{2k_{0}}}.
\end{align}
Similarly, we can get as $k_{0}\rightarrow \infty$,
\begin{align} \label{esbe1}
	\beta_{1}(\tilde{f})=\int_{\Delta_{2k_{0}}(\frac{2}{3\theta_{0}})}t_{2}^{2} \left(\frac{\partial^{2k_{0}+1}\tilde{f}(t_{1},\cdot\cdot\cdot,t_{2k_{0}})}{\partial t_{1}(\partial t_{2})^{2}\cdot\cdot\cdot \partial t_{2k_{0}}} \right)^{2}\:\mathrm{d} \underline{t}
	=o\left( \frac{\gamma^{2k_{0}-2}}{(2k_{0})^{2k_{0}}}\right)
\end{align}
and 
\begin{align} \label{presbe2}
	\int_{\Delta_{2k_{0}}(\frac{2}{3\theta_{0}})}t_{2}\left(\frac{\partial^{2k_{0}}\tilde{f}(t_{1},\cdot\cdot\cdot,t_{2k_{0}})}{\partial t_{1}\partial t_{2}\cdot\cdot\cdot \partial t_{2k_{0}}} \right)^{2}\:\mathrm{d}\underline{t}
	=o\left( \frac{\gamma^{2k_{0}-2}}{(2k_{0})^{2k_{0}}}\right).
\end{align}
By Cauchy-Schwarz we deduce from \eqref{esarp} and \eqref{presbe2} that
\begin{align} \label{esbe2}
	\beta_{2}(\tilde{f})=\int_{\Delta_{2k_{0}}(\frac{2}{3\theta_{0}})}t_{2}\left(\frac{\partial^{2k_{0}+1}\tilde{f}(t_{1},\cdot\cdot\cdot,t_{2k_{0}})}{\partial t_{1}(\partial t_{2})^{2}\cdot\cdot\cdot \partial t_{2k_{0}}} \right)\left(\frac{\partial^{2k_{0}}\tilde{f}(t_{1},\cdot\cdot\cdot,t_{2k_{0}})}{\partial t_{1}\partial t_{2}\cdot\cdot\cdot \partial t_{2k_{0}}} \right)\:\mathrm{d}\underline{t}
	=o\left(\frac{\gamma^{2k_{0}-2}}{(2k_{0})^{2k_{0}}}\right),
\end{align}
as $k_{0}\rightarrow \infty$. We conclude that from \eqref{ATT}, \eqref{esarp}, \eqref{esbe1} and \eqref{esbe2},
 \begin{align} \label{s2conpds}
	S_{2}
	\leq &\frac{k_{0}x\log x}{(\log R)^{2k_{0}+1}}\frac{W^{2k_{0}-1}}{\phi(W)^{2k_{0}}}\left(\frac{0.168}{\frac{2}{3\theta_{0}}-1}\cdot\frac{\gamma^{2k_{0}-2}}{(2k_{0})^{2k_{0}}}\right),
\end{align}
provided $k_{0}$ and $x$ are large.

\begin{remark}
	In \cite[eq. (3.19)]{L-P2015}, Li and Pan established $\alpha\leq 8.98 \left(\frac{\gamma}{2k_{0}}\right)^{2k_{0}-1}$. In comparison to their result, the order of our upper bound for $\alpha$ (cf. \eqref{esarp}) saves a factor $k_{0} \gamma \sim k_{0}/\log k_{0}$. This distinction highlights the power of Tao's approach and plays a crucial role in deriving our main theorem.
\end{remark}
\subsection{Completion of the proof of Theorem \ref{NewTh}} \label{ComP}

\

Plugging the estimates for $S_{1}$, $S_{2}$, and $S_{3}$ (see \eqref{s1finaLBPRS}, \eqref{s2conpds}, and \eqref{s3FinaLBNP}) into  \eqref{rewS} yields
\begin{align*} 
	S(x, C) \geq& \frac{\gamma^{2k_{0}}}{(2k_{0})^{2k_{0}}}\frac{x}{(\log R)^{2k_{0}-1} \log x} \frac{W^{2k_{0}-1}}{\phi(W)^{2k_{0}}}
	\\
	&\cdot\left(\frac{2m}{\theta_{0}} +1.064(1-c_{1})-\frac{k_{0}}{C \gamma^{2} }\left(\ \frac{\log x}{\log R}\right)^{2} \frac{0.168}{\frac{2}{3\theta_{0}}-1}-m \frac{\log x}{ \log R} +o(1)\right),
\end{align*}
as $x \rightarrow \infty$.
Since $R=x^{\frac{\theta_{0}}{2}-\frac{1}{100000m}}$, $\theta_{0}>0.5$, $c_{1}<8\times 10^{-6}$, and $m$ is large, we have 
\[\frac{2m}{\theta_{0}} +1.064(1-c_{1})-m \frac{\log x}{ \log R} > 1.063(1-c_{1}).
\]
We therefore get $S(x, C)>0$ for large $x$ when 
\begin{align} \label{chosC}
	C= \frac{k_{0}}{\gamma^{2}}\cdot \frac{0.168}{(\frac{\theta_{0}}{2}-\frac{1}{100000m})^{2}(\frac{2}{3\theta_{0}}-1)}\cdot\frac{1}{1.063(1-c_{1})}.
\end{align}
It follows from \eqref{Chk0} that
\[\log k_{0}=2\log m+\frac{4m}{\theta_{0}(1-c_{1})}+8.
\]
Recall $\gamma=\frac{1}{A}(1-\frac{1}{1+AT})$, $A=\log 2k_{0}-2\log \log 2k_{0}$, and $T=(e^{A}-1)/A$. We find that
\begin{align*}
	&-\log \gamma 
	\\
	= &\log \log k_{0} + \log \left(1+ \frac{\log 2}{\log k_{0}}\right) + \log \left(1- \frac{2 \log \log 2k_{0}}{\log 2k_{0}}\right)-\log \left(1-\frac{1}{1+AT}\right) 
	\\
	= & \log m + \log \left( \frac{4}{\theta_{0}(1-c_{1})}+ \frac{2 \log m +8}{m}\right) + \log \left(1+ \frac{\log 2}{\log k_{0}}\right)+ \log \left(1- \frac{2 \log \log 2k_{0}}{\log 2k_{0}}\right)
	\\
	&\quad -\log \left(1-\frac{1}{1+AT}\right) = \log m + \log \left(\frac{4}{\theta_{0}(1-c_{1})}\right)+o(1),
	 \ \mbox{as} \ k_{0} \rightarrow \infty.
\end{align*}
Hence,
\begin{align} \label{fiea1}
	\log k_{0} -2\log \gamma &= 2 \log m + \frac{4m}{\theta_{0}(1-c_{1})}+8 + 2\log m + 2\log \left(\frac{4}{\theta_{0}(1-c_{1})}\right)+ o(1)
	\notag
	\\
	&= 4 \log m + \frac{4m}{\theta_{0}(1-c_{1})}+8 + 2\log \left(\frac{4}{\theta_{0}(1-c_{1})}\right)+o(1),
	\ \mbox{as} \ k_{0} \rightarrow \infty.
\end{align}
Combining \eqref{chosC} and \eqref{fiea1} gives that $\log C$ is bounded by
\begin{align*}  
&\log k_{0}- 2\log \gamma + \log \frac{0.168}{1.063(1-c_{1})(\frac{\theta_{0}}{2}-\frac{1}{100000m})^{2}(\frac{2}{3\theta_{0}}-1)}
	\\
	<& 4 \log m + \frac{4m}{\theta_{0}(1-c_{1})}+8 + 2\log \left(\frac{4}{\theta_{0}(1-c_{1})}\right)+  \log \frac{0.168}{1.063(1-\frac{8}{ 10^{6}})\cdot\frac{1}{16}\cdot(\frac{2\times 1318}{3 \times 691}-1)}
	\\
	\leq &7.63m + 4 \log m + 21 \log 2
\end{align*}
for large $k_{0}$. Here we used $\frac{4}{\theta_{0}(1-c_{1})} < 7.63$ and $\frac{1}{2} < \theta_{0} < \frac{691}{1318}$ (cf. \eqref{snres}). The proof of Theorem \ref{NewTh} is now complete in view of the discussion in the Section \ref{SU}.

\section{Further remarks}
In Stadlmann's recent work \cite{JuS}, she obtained a new minorant for the indicator function of primes, achieving a sharper exponent of distribution to smooth moduli, 0.5253, compared to Baker and Irving's construction. By substituting this new minorant (see \cite[Proposition 2 on page 15]{JuS}) for Lemma \ref{BIM}, one can derive an improved upper bound of $q_{m+1}-q_{m}=O(e^{7.615m})$ in Theorem \ref{NewTh}.

Due to Chen's celebrated theorem \cite{Chen-1978}, we call $p$ a Chen prime if $\Omega(p+2)\leq 2$.
With current methods, establishing bounded gaps between Chen primes still appears to be out of reach.

\section*{Acknowledgements}
The author would like to thank to thank Terence Tao and Gregory Debruyne for many useful conversations and suggestions. The work was supported by the China Scholarship Council (CSC).

\end{document}